\documentclass{article}
\usepackage{amsmath,amssymb,amsfonts}

\newtheorem{lemma}{Lemma}
\newtheorem{theorem}{Theorem}

\newtheorem{corollary}{Corollary}

\def\D{{\cal D}}
\def\C{{\mathbb C}}
\def\R{{\mathbb R}}
\def\Z{{\mathbb Z}}
\def\M{{\cal M}}

\def\Im{{\mathrm{Im}\, }}

\def\Ker{{\mathrm{Ker}\,}}

\def\beq{\begin{equation}}
\def\eeq{\end{equation}}

\begin{document}

\title{Infinitesimal Darboux transformations of
the spectral curves of tori in the four-space \thanks{The work was
supported by RFBR (grants no. 05-01-01032a (P.G.G.) and 06-01-00094a
(I.A.T.)), by the program ``Fundamental problems of nonlinear
dynamics'' of the Presidium of RAS, and by the program ``Leading
scientific schools'' (grant NS-4182.2006.1). The second author was
also supported by the complex integration project 2.15 of SB RAS.}}
\author{P.G. Grinevich
\thanks{Landau Institute of Theoretical Physics, Kosygin street 2,
117940 Moscow, Russia; e-mail: pgg@landau.ac.ru.} \and I.A.
Taimanov
\thanks{Institute of Mathematics, 630090 Novosibirsk, Russia;
e-mail: taimanov@math.nsc.ru}}
\date{}

\maketitle

\section{Introduction}
\label{sec1}

The present paper resumes the study of relations between the
geometric properties of surfaces in $\R^3$ and $\R^4$ and the
spectral properties of the corresponding Dirac operators started in
\cite{T1}. In the paper we study the behavior of the spectral curve
of a torus in $\R^4$ under conformal transformations of $\bar{\R}^4$
and, in particular, prove that

\begin{itemize}
\item
the conformal transformations of $\bar{\R}^4$ which map a torus $T
\subset \R^4$ into a compact torus preserve all Floquet multipliers
of the corresponding Dirac operator
\end{itemize}

\noindent (see Theorem 1 and Corollary 1).

This generalizes the analogous result for tori in $\R^3$  which
was proved by M.U. Schmidt and the first author (P.G.G.) \cite{GS}
and confirmed the conjecture of the second author (I.A.T.) on the
conformal invariance of the spectral curves of tori in $\R^3$.
Therewith by the spectral curve it was understood the analytic set
$\M (\Gamma)$ in $\C^2$ formed by the Floquet multipliers  on the
zero energy level.

There is a more careful and precise definition of the spectral
curve compatible with the finite gap integration theory. It reads
that this curve $\Gamma$ gives a one-to-one parameterization of a
basis for Floquet--Bloch functions of the corresponding Dirac
operator and the multiplier set $\M (\Gamma)$ is the image of the
mapping $\M$ which corresponds to each point the multipliers of
the corresponding Floquet function (see \cite{T3}). Therewith we
show that

\begin{itemize}
\item
the spectral curve of a Dirac operator corresponding to a torus $T
\subset \R^4$ is not always preserved by conformal transformations
of $\R^4$. However possible deformations consist in gluing together
or ungluing points with some fixed multipliers $(\kappa_1,\kappa_2)$
\end{itemize}

\noindent (see Remark 1).

This almost isospectrality effect was not noticed before and we
think that it is interesting by itself. It does not relate to the
dimension of the ambient space and holds also for tori in the
three-space. In fact the nontrivial actions of conformal
transformations of $\R^3$ or $\R^4$ on the potential of the Dirac
operator corresponding to an immersed torus are described by
nonlinear systems of the Melnikov type. Thus

\begin{itemize}
\item
the Melnikov type deformations may be only almost isospectral and we
show that in some important geometrical examples (related to the
Clifford tori) after reparameterization  of the temporary variable
such flows reduce to integrable systems on whiskered tori
\end{itemize}

\noindent (see \S 5.1).

To the spectral curve of the Dirac operator with given local
parameters near ``infinities'' there corresponds an infinite
family of nonlocal conservation laws which for the case of the
one-dimensional potential reduce to the Kruskal--Miura integral.
Since their values in fact depend on $\M(\Gamma)$ these integrals
are also preserved by the infinitesimal conformal deformations.

\section{The Weierstrass representation and the spectral curve}

\subsection{The Weierstrass representation}

The Weierstrass representation of surfaces in $\R^3$ and $\R^4$
there corresponds to a surface the Dirac operator with potentials
\beq \label{dirac} \D = \left(
\begin{array}{cc}
0 & \partial \\
-\bar{\partial} & 0
\end{array}
\right) + \left(
\begin{array}{cc}
U & 0 \\
0 & \bar{U}
\end{array}
\right) \eeq where $U = U(z,\bar{z})$ and $z$ is a conformal
parameter on surface. In the sequel we shall use the following
agreement: we write $f(z)$ instead of $f(z,\bar z)$ and the
notation $f(z)$ does not imply that $f$ is holomorphic unless the
opposite is not stated explicitly.

In fact, $U$ is defined as a section of some bundle over the surface.
In the paper we consider the case of tori.
By the uniformization theorem, every torus is conformally equivalent
to a flat torus $\R^2/\Lambda$
where $\Lambda$ is the period lattice.

We have
\begin{itemize}
\item
for tori in $\R^3$ the potential $U$ is double-periodic:
$U(z+\gamma) = U(z)$ for all $\gamma \in \Lambda$, and real-valued;

\item
for tori in $\R^4$ the potential $U$ is double-periodic
however is defined up to transformations
\beq
\label{curvetransform}
U \to U e^{\bar{a} + \overline{bz} - a -bz}
\eeq
where $a,b \in \C$ and $\Im b\gamma \in \pi \Z$ for all $\gamma \in
\Lambda$.
\end{itemize}

A surface in $\R^4$ is given by the integral formulas
$$
x^k = x^k_0 + \int \left( x^k_z dz + \bar{x}^k_z
d\bar{z}\right), \ \ k=1,2,3,4,
$$
with
$$
x^1_z = \frac{i}{2} (\bar{\varphi}_2\bar{\psi}_2 + \varphi_1
\psi_1), \ \ \ \ x^2_z = \frac{1}{2} (\bar{\varphi}_2\bar{\psi}_2 -
\varphi_1 \psi_1),
$$
$$
x^3_z = \frac{1}{2} (\bar{\varphi}_2 \psi_1 + \varphi_1
\bar{\psi}_2), \ \ \ \ x^4_z = \frac{i}{2} (\bar{\varphi}_2 \psi_1 -
\varphi_1 \bar{\psi}_2)
$$
where $x^1,\dots,x^4$ are the Euclidean coordinates in $\R^4$ (see
\cite{K2,PP}) and $\psi$ and $\varphi$ meet the Dirac equations
$$
\D \psi = \left( \begin{array}{cc}  U & \partial \\ -\bar\partial &
\bar U  \end{array}  \right)\psi = 0, \ \ \
\D^\vee \varphi = \left( \begin{array}{cc}  \bar U & \partial \\
-\bar\partial & U  \end{array}  \right) \varphi = 0
$$
(we note that the operators $\D$ and $\D^\vee$ are Hermitian
conjugate).

The case of surfaces in $\R^3$ is obtained as the reduction:
$U=V=\bar{U}$, $\D = \D^\vee$, and for $\psi = \varphi$ we obtain
a surface lying in the linear subspace $x^4=0$.

Although for surfaces in $\R^3$
the derivation of the Weierstrass representation for surfaces in $\R^3$
is straightforward
and the vector function $\psi=\varphi$ is easily defined from the
geometrical data \cite{T1}
for surfaces in $\R^4$ the situation is different.
This procedure is much more delicate and the functions $\psi$,
$\varphi$ and
$U$ are defined from nonlinear equations which have to be solved
globally on the whole surface \cite{T2}.

\subsection{The spectral curve}
\label{sec2}

In the framework of differential geometry of surfaces the spectral
curves appear as the spectral curves
of integrable surfaces (constant mean curvature tori in $\R^3$
\cite{Bobenko} and harmonic tori
in $S^3$ \cite{Hitchin}) and were used for constructing the
explicit formulas for such tori in terms of theta functions
corresponding to
these spectral curves which appear to be of finite genus.
Recently they were used for obtaining the lower estimates for the areas
of minimal tori in $S^3$
\cite{FLPP}.

For general tori in $\R^3$ the spectral curve was defined via the
Weierstrass representation by the second author as the spectral
curve of the Dirac operator associated, i.e. with the potential $U =
e^\alpha H/2$ where $e^{2\alpha}dz d\bar{z}$ is the induced metric
and $H$ is the mean curvature.

First the spectral curve on the zero energy level was introduced by
Dubrovin, Krichever, and Novikov for the two-dimensional
Schr\"odinger operator \cite{DKN}.

The general definition of the Floquet--Bloch spectrum of a
multi-periodic operator was introduced in \cite{Novikov}. For the
two-dimensional case it is as follows.
Given a double-periodic operator $L$, its Floquet\,(--Bloch) function
$\psi$
is defined as the formal solution to the equation
$$
L \psi = E\psi
$$
meeting the periodicity conditions
$$
\psi(z + \gamma_i) = \kappa_i \psi(z), \ \ \ i=1,2,
$$
where $\gamma_1,\gamma_2$ are the generators of the period lattice.
It is said that $E$ is the energy level and $\kappa_1,\kappa_2 \in
\C \setminus \{0\}$ are the (Floquet) multipliers of $\psi$. Let us
put $E=0$, i.e., let us consider the zero energy level, and assume
that the possible values of the multipliers meet some analytical
dispersion relation
$$
F(\kappa_1,\kappa_2)=0
$$
which defines a one-dimensional complex manifold $\M (\Gamma)$
(complex curve).

The dispersion relations do exist not for all operators. However
this picture is true for elliptic operators and some other operators
closed to them. In the middle of 1980s the problem of rigorous
confirmation of this picture, in particular, for the two-dimensional
Schr\"o\-din\-ger operator and the heat operator was addressed by
two different ways:

1) by perturbation methods Krichever did construct the spectral
curve and the Floquet--Bloch eigenfunction as perturbations of their
counterparts for the operator with the zero potential
\cite{Krichever}. Therewith for the case of the Sch\"rodinger
operator we have two infinite ends at which the eigenfunction
asymptotically behaves as a holomorphic function at one end and as a
antiholomorphic function at another end and the perturbation of the
spectral curve consists in opening resonant pairs into handles
outside some compact part at which the perturbation may result in
more complicated topological surgery. This geometrical picture
rising to the spectral theory initiated the development of the
analytical theory of such Riemann surfaces (non only hyperelliptic)
of infinite genus \cite{FKT}.

2) the second author (I.A.T.) demonstrated how to obtain the
analytical dispersion relation for hyperelliptic periodic operators
by using the Fredholm alternative for analytical pencils of
operators (the Keldysh theorem, for the Dirac operator such a proof
is exposed in \cite{T3}). Therewith $F(\kappa_1,\kappa_2)$ is a
(regularized) determinant of the operator $L(\kappa_1,\kappa_2)$
which is the operator $L$ defined on the space of functions meeting
the boundary conditions \beq \label{bloch-1}
\psi(z+\gamma_1)=\kappa_1 \psi(z), \ \ \psi(z+\gamma_2)=\kappa_2
\psi(z), \ \ \psi=\left(
\begin{array}{c} \psi_1 \\ \psi_2
\end{array} \right).
\eeq

 In general to several points of this manifold there
correspond not a one-dimensional family of Floquet eigenfunctions
with these multipliers. Therefore to obtain the spectral curve
$\Gamma$ we have to consider a partial normalization of $\M
(\Gamma)$. If it is of finite genus it admits a compactification by
finitely many points (the ``infinities'') to an algebraic curve.

Therefore, we have to make difference between

\begin{itemize}
\item
the spectral curve $\Gamma$ such that there is a rank one bundle,
over $\Gamma$, which is formed by the corresponding Floquet
functions such that any Floquet function is a linear composition
of the generators of fibres \footnote{This is exactly the complex
curve which is coming into the finite gap integration scheme.};

\item
the multiplier set $\M = \M(\Gamma)$ which parameterizes the pairs of
Floquet multipliers $(\kappa_1,\kappa_2)$.
This set is the image of the multiplier mapping
$$
\M: \Gamma \to \C^2, \ \ \ \M(\psi) = (\kappa_1,\kappa_2);
$$

\item
the algebraic spectral curve which is obtained by a compactification by
finitely many points from
the curve $\Gamma$ of finite genus.
\end{itemize}

The difference between these three spectral curves is demonstrated
in \S \ref{subsec5.2} .for the one-dimensional Schr\"odinger
operator with a complex-valued potential.

We refer for the detailed explanation of this picture for the
Dirac operator to \cite{T3}. This survey also explains the
approach which was proposed by the second author (I.A.T.) to
proving the Willmore conjecture by using the spectral curves. In
part, this approach comes from the observation that the Willmore
functional
$$
{\cal W}(M) = \int_M H^2 d\mu
$$
for surfaces $M$ in $\R^3$ is related to their
Weierstrass representation and, in particular, to the spectral curve
(here $H$ is the mean
curvature and $d\mu$ is the induced measure on the surface).

For closed surfaces the Willmore functional is invariant under
conformal transformations of $\bar{\R}^3$ which preserve the
compactness of the surface. This led to another conjecture that
the spectral curve for tori is invariant with respect to these
transformations and this conjecture was confirmed in \cite{GS}. On
the modern language we may say that in \cite{GS} it was proved
that the multiplier set $\M$ is invariant.

Hence, given the potential $U$ for a torus in $\R^4$, we have a pair
of spectral curves $\Gamma$ and $\Gamma^\vee$, i.e., the spectral
curves of the operators $\D$ and $\D^\vee$. These spectral curves
are closely related and, in particular, we have the evident

\begin{lemma}
The operators $\D(\kappa_1,\kappa_2)$ and
$\D^\vee(\bar\kappa_1,\bar\kappa_2)$ are Hermitian conjugate for all
$\kappa_1$, $\kappa_2$ and their indices are equal to $0$.

Therefore

1) $\dim \Ker \D(\kappa_1,\kappa_2) = \dim \Ker
\D^\vee(\bar\kappa_1,\bar\kappa_2)$ for all $\kappa_1, \kappa_2$;

2) The multiplier sets for the operators $\D$ and $\D^\vee$ are
complex conjugate.
\end{lemma}

We also have

\begin{lemma}
The spectral curves for  operators $\D$ and $\D^\vee$ are both real.
\end{lemma}

{\it Proof}. Let
$$
\psi=\left( \begin{array}{c}  \psi_1  \\ \psi_2  \end{array} \right)
$$
be a zero Floquet eigenfunction of $\D$ with multipliers  $\kappa_1$,
$\kappa_2$. Then
the function
$$
\psi^\ast = \left( \begin{array}{c} \bar \psi_2  \\ -\bar\psi_1
\end{array} \right)
$$
is also a zero Floquet eigenfunction of $\D$ with multipliers
$\bar\kappa_1$, $\bar\kappa_2$. Therefore $\Gamma$ admits the
antiholomorphic involution $\psi \to \psi^\ast, (\kappa_1,\kappa_2)
\rightarrow
(\bar\kappa_1,\bar\kappa_2)$.

We remark that the potential $U$ of a torus in $\R^4$ is defined
up to transformations (\ref{curvetransform}) which imply the
following simple transformations of the Floquet--Bloch
eigenfunctions $\psi$ of $\D$ and their multipliers:
$$
\left(\begin{array}{c}
\psi_1 \\ \psi_2 \end{array}\right) \to
\left(\begin{array}{c}
e^{a+bz} \psi_1 \\ e^{\bar{z} + \bar{b}\bar{z}}\psi_2
\end{array}\right), \ \ \
\kappa_i \to e^{b\gamma_i}\kappa_i, \ \ \ i=1,2
$$
(see \cite{T2}).

Therefore the Dirac operator $\D$ corresponding to a conformal
immersion of a torus $\R^2/\Lambda$ into $\R^4$ is not unique: we
may put for it $\D$ or $\D^\vee$ or any of their transformations of
the form (\ref{curvetransform}). If the torus lies in $\R^3 \subset
\R^4$ there is a normalization condition $U = \bar{U}$ which fixes
the operator uniquely. However transformations of the form
(\ref{curvetransform}) as well as the changes of the base for the
period lattice $\Lambda$ does not change the spectral curve and only
transforms the multiplier mapping $\M$.

\section{The infinitesimal Darboux transformation}
\label{sec3}

Let
\begin{equation}
\label{operators-1}
 L=\left( \begin{array}{cc} \bar \partial  &  -p \\  -q &  \partial
\end{array}  \right), \ \
 L^\vee=\left( \begin{array}{cc}  -\bar\partial  & - q \\  -p &
\partial \end{array}  \right), \ \
p=p(z), \ \ q=q(z).
\end{equation}
be a pair of linear operators with periodic coefficients:
$$
p(z+\gamma_1)=p(z+\gamma_2)=p(z), \ \
q(z+\gamma_1)=q(z+\gamma_2)=q(z),
$$
i.e. the functions $p$, $q$ are defined on a torus $T$  with the
periods  $\gamma_1$ and $\gamma_2$.

Let us define {\sl the infinitesimal Darboux transformation}.

Let us consider

\begin{itemize}
\item
$\Psi^D$ and $\Phi^D$ be a pair of zero Floquet eigenfunctions of
$L$ and $L^\vee$ respectively such that
$$
L\Psi^D=0, \ \  L^\vee \Phi^D =0, \ \
$$
$$
\Psi^D(z+\gamma_i)=\hat\kappa_i \Psi^D(z), \ \
\Phi^D(z+\gamma_i)=1/\hat\kappa_i \Phi^D(z), \ \ i=1,2,
$$
$$
\Psi^D=\left( \begin{array}{c}  \Psi^D_1(z)  \\ \Psi^D_2(z)
\end{array} \right), \\
\Phi^D=\left( \begin{array}{c}  \Phi^D_1(z)  \\ \Phi^D_2(z)
\end{array} \right);
$$

\item
a family of  zero Floquet--Bloch eigenfunctions  $\psi(\lambda,z)$
of $L$ with multipliers $\kappa_1(\lambda)$, $\kappa_2(\lambda)$
respectively, $\lambda\in\Gamma$;

\item
a family of zero Floquet--Bloch eigenfunctions of $L^\vee$:
$\phi(\mu,z)$ with
multipliers $\kappa^\vee_1(\mu)$, $\kappa^\vee_2(\mu)$
respectively, $\mu\in\Gamma^\vee$.
\end{itemize}

Let us define the following pair of forms:

\begin{equation}
\label{omega1} d\omega(\lambda,z) =\Phi^D_1(z)
\psi_1(\lambda,z)dz-\Phi^D_2(z)\psi_2(\lambda,z)d\bar z
\end{equation}
\begin{equation}
\label{omega2} d\omega^\vee(\mu,z) =
\phi_1(\mu,z)\Psi^D_1(z)dz-\phi_2(\mu,z)\Psi^D_2(z)d\bar z
\end{equation}

\begin{lemma}
The both forms $d\omega(\lambda,z)$ and $ d\omega^\vee(\mu,z)$ are
closed. Therefore formulas (\ref{omega1}) and (\ref{omega2})
defines the functions  $\omega(\lambda,z)$ and
$\omega^\vee(\mu,z)$ up to integration constants $c(\lambda)$ and
$c^\vee(\mu)$.
\end{lemma}

{\it Proof} is straightforward:
$$
d(d\omega)=-(\partial_{\bar z}[\Phi^D_1 \psi_1(\lambda)] +
\partial_{z}[(\Phi^D_2 \psi_2(\lambda)]) dz\wedge d\bar z=
$$
$$
-( p\Phi^D_2 \psi_1(\lambda)-q\Phi^D_1 \psi_2(\lambda) + q
\Phi^D_1 \psi_2(\lambda) - p\Phi^D_2 \psi_1(\lambda)  )dz\wedge
d\bar z =0.
$$

\begin{lemma}
\label{Bloch-Cauchy} Let either $\kappa_1(\lambda)/\hat\kappa_1\ne
1$ or $\kappa_2(\lambda)/\hat\kappa_2\ne 1$. Then the function
$\omega(\lambda,z)$ is defined uniquely by the Floquet--Bloch
condition:
$$
\omega(\lambda,z+\gamma_i)=\frac{\kappa_i(\lambda)}{\hat\kappa_i}
\omega(\lambda,z).
$$
If either $\kappa^\vee_1(\mu) \hat\kappa_1\ne 1$ or
$\kappa^\vee_2(\mu) \hat\kappa_2\ne 1$, the function
$\omega^\vee(\mu,z)$ is defined uniquely by the Floquet--Bloch
condition:
$$
\omega^\vee(\mu,z+\gamma_i)=\kappa^\vee_i(\mu)\hat\kappa_i
\omega^\vee(\mu,z).
$$
\end{lemma}

{\it Proof}. Let us consider the basic parallelogram $0, \gamma_1,
\gamma_1+\gamma_2, \gamma_2$. Let us denote
$$
I_1(\lambda)=\int\limits_{0}^{\gamma_1}d\omega(\lambda), \ \
I_2(\lambda)=\int\limits_{0}^{\gamma_2}d\omega(\lambda).
$$
We have
$$
\int\limits_{\gamma_2}^{\gamma_1+\gamma_2}d\omega(\lambda) =
I_1(\lambda) \cdot \kappa_2(\lambda)/\hat\kappa_2, \ \ \ \
\int\limits_{\gamma_1}^{\gamma_1+\gamma_2}d\omega(\lambda) =
I_2(\lambda) \cdot \kappa_1(\lambda)/\hat\kappa_1.
$$
Therefore
$$
I_1(\lambda)+I_2(\lambda) \cdot \kappa_1(\lambda)/\hat\kappa_1 =
I_2(\lambda) + I_1(\lambda) \cdot \kappa_2(\lambda)/\hat\kappa_2,
$$
and
\begin{equation}
\label{compat-bloch} I_1(\lambda)[ \kappa_2(\lambda)/\hat\kappa_2
- 1 ] = I_2(\lambda) [ \kappa_1(\lambda)/\hat\kappa_1 - 1].
\end{equation}
The periodicity condition for $\omega(\lambda)$ implies
$$
\omega(\lambda,\gamma_1)=\omega(\lambda,0)\kappa_1(\lambda)/\hat\kappa_1
= \omega(\lambda,0) + I_1(\lambda),
$$
$$
\omega(\lambda,\gamma_2)=\omega(\lambda,0)\kappa_2(\lambda)/\hat\kappa_2
= \omega(\lambda,0) + I_2(\lambda),
$$
\begin{equation}
\label{omega-bloch}
\begin{split}
\omega(\lambda,0)\cdot [ \kappa_1(\lambda)/\hat\kappa_1 - 1 ] =
I_1(\lambda), \\
\omega(\lambda,0)\cdot [ \kappa_2(\lambda)/\hat\kappa_2 - 1 ] =
I_2(\lambda).
\end{split}
\end{equation}
From (\ref{compat-bloch}) it follows, that under the assumptions
of Lemma the system (\ref{omega-bloch}) is compatible and has an
unique solution.

\begin{theorem}
\label{lemma-iso}
Consider the following variation of the space of
Floquet--Bloch functions:
\begin{equation}
\label{bloch-wave} \left\{
\begin{array}{l}
\delta\psi_1(\lambda,z)=\omega(\lambda,z)\Psi^D_1(z)\\
\delta\psi_2(\lambda,z)=\omega(\lambda,z)\Psi^D_2(z)\\
\delta\Phi^D_1(\mu,z)=\omega^\vee(\mu,z)\Phi^D_1(z)\\
\delta\Phi^D_2(\mu,z)=\omega^\vee(\mu,z)\Phi^D_2(z)
\end{array}\right.
\end{equation}
\begin{enumerate}
\item This deformation corresponds to the following variation of the
operators $L$, $L^\vee$:
\begin{equation}
\label{melnikov-1}
\begin{split}
\delta L =\left( \begin{array}{cc} 0  & \Psi^D_1(z)\Phi^D_2(z) \\
-\Psi^D_2(z)\Phi^D_1(z)  & 0
\end{array} \right) , \\
\delta L^\vee=\left( \begin{array}{cc} 0
& -\Psi^D_2(z)\Phi^D_1(z)
\\  \Psi^D_1(z)\Phi^D_2(z) & 0 \end{array}  \right), \ \
\end{split}
\end{equation}
Therefore it is self-consistent (variations of different wave
functions result in the same variation of potentials), and
respects the symmetry between $L$, $L^\vee$. In terms of $p$, $q$
we get Melnikov-type variations of potentials \cite{M}:
\begin{equation}
\delta p(z) = - \Psi^D_1(z)\Phi^D_2(z), \ \ \delta q(z) =
\Psi^D_2(z)\Phi^D_1(z).
\end{equation}

\item For all $\lambda$, $\mu$ such, that the conditions of
Lemma~\ref{Bloch-Cauchy}
are fulfilled it is natural to normalize the kernels
$\omega(\lambda)$, $\omega^\vee(\mu)$ by the Floquet--Bloch
conditions. Then deformation (\ref{bloch-wave}) respects the
Floquet--Bloch properties of the functions $\psi(\lambda)$,
$\phi(\mu)$, and does not change the multipliers.
\end{enumerate}
\end{theorem}

{\it Proof}. From (\ref{operators-1}) it follows:
$$
p=\partial_{\bar z}\psi_1(\lambda)/\psi_2(\lambda), \ \
q=\partial_{z}\psi_2(\lambda)/\psi_2(\lambda),
$$
$$
\delta p=\frac{\partial_{\bar
z}[\delta\psi_1(\lambda)]}{\psi_2(\lambda)}- \frac{[\partial_{\bar
z}\psi_1(\lambda)]\delta\psi_2(\lambda)}{\psi^2_2(\lambda)}=
$$
$$
=\frac{\partial_{\bar
z}[\omega(\lambda)\Psi^D_1]}{\psi_2(\lambda)} - p
\frac{\omega(\lambda)\Psi^D_2}{\psi_2(\lambda)}=
\frac{[\partial_{\bar
z}\omega(\lambda)]\Psi^D_1}{\psi_2(\lambda)}+
\frac{\omega(\lambda)\partial_{\bar z}\Psi^D_1}{\psi_2(\lambda)} -
p \frac{\omega(\lambda)\Psi^D_2}{\psi_2(\lambda)}=
$$
$$
=\frac{-\psi_2(\lambda)\Phi^D_2\Psi^D_1}{\psi_2(\lambda)}+
\frac{\omega(\lambda)p\Psi^D_2}{\psi_2(\lambda)} - p
\frac{\omega(\lambda)\Psi^D_2}{\psi_2(\lambda)}= -
\Psi^D_1\Phi^D_2
$$
$$
\delta q=\frac{\partial_{z}[\delta\psi_2(\mu)]}{\psi_1(\mu)}-
\frac{[\partial_{z}\psi_2(\mu)]\delta\psi_1(\mu)}{\psi^2_1(\mu)}=
$$
$$
=\frac{\partial_{z}[\omega^\vee(\mu)\Psi^D_2]}{\psi_1(\mu)} - q
\frac{\omega^\vee(\mu)\Psi^D_1}{\psi_1(\mu)}=
\frac{[\partial_{z}\omega^\vee(\mu)]\Psi^D_2}{\psi_1(\mu)}+
\frac{\omega^\vee(\mu)\partial_{z}\Psi^D_2}{\psi_1(\mu)} - q
\frac{\omega^\vee(\mu)\Psi^D_1}{\psi_1(\mu)}=
$$
$$
=\frac{\psi_1(\mu)\Phi^D_1\Psi^D_2}{\psi_1(\mu)}+
\frac{\omega^\vee(\mu)q\Psi^D_1}{\psi_1(\mu)} - q
\frac{\omega^\vee(\mu)\Psi^D_1}{\psi_1(\mu)}= \Psi^D_2\Phi^D_1.
$$

Theorem \ref{lemma-iso} and, in particular, part 2 shows that the
deformation respects the Floquet--Bloch properties of the function
$\psi(\lambda,z)$ if the condition of Lemma~\ref{Bloch-Cauchy} is
fulfilled which, in particular, means that
$\kappa_1(\lambda)/\hat{\kappa}_1 \neq 1 \ \ \mbox{or} \ \ \
\kappa_2(\lambda)/\hat{\kappa}_2 \neq 1.$ These inequalities are
fulfilled at a generic point however there is a discrete set of
points of the spectral curve at which we have \beq \label{Lset}
 \kappa_1(\lambda)/\hat{\kappa}_1 =
\kappa_2(\lambda)/\hat{\kappa}_2 = 1. \eeq Hence this deformation
it preserves the multiplier set of $L$ outside of a discrete set
of points satisfying (\ref{Lset}). Since the multiplier set is
analytic, it is preserved. Of course, the same becomes valid for
$L^\vee$ and its Floquet--Bloch functions $\phi(\mu,z)$ after
replacing (\ref{Lset}) by the condition \beq \label{Lveeset}
\kappa^\vee_1(\mu)\hat{\kappa}_1 =
\kappa^\vee_2(\mu)\hat{\kappa}_2 = 1. \eeq

\begin{corollary}
The infinitesimal deformation (\ref{melnikov-1}) preserves the
multiplier sets $\M (\Gamma)$ (on the zero energy level) for the
operators $L$ and $L^\vee$.
\end{corollary}

{\sc Remark 1.} Theorem \ref{lemma-iso} does not imply that the
spectral curve is preserved. Indeed the forms $\omega$ and
$\omega^\vee$ are not defined at points meeting the condition
(\ref{Lset}). Given a  multiple point meeting this condition the
analytic continuation of $\omega$ and $\omega^\vee$ a priori gives
its own limit at each branch. Therewith the deformation is
correctly defined on the normalization of $\Gamma$ and the
corresponding Floquet function evolves differently and this leads
to decreasing the multiplicity of a singular point on the spectral
curve. Of course, the converse is also possible. The examples from
\S \ref{sec5} shows that that may take place and how that happens.

{\sc Remark 2.} The idea of using kernels analogous to
(\ref{omega1}) and (\ref{omega2}) for calculating the
Floquet--Bloch functions deformations for all values of spectral
parameter was first suggested in \cite{GO} by A.Yu. Orlov and the
first author (P.G.G.).

\section{Proof of the conformal invariance}
\label{sec4}

Let us apply Theorem \ref{Bloch-Cauchy} to conformal
transformations of tori in $\R^4$ induced by conformal
transformations of the ambient space $\bar{\R}^4$.

We assume  we have a conformal immersion of torus $\R^2/\Lambda$ into
${\mathbb R}^4$
defined by:
\begin{equation}
\label{weier4}
\begin{split}
\partial_z x^1 = \frac{i}{2} (\bar{\Phi}_2\bar{\Psi}_2 + \Phi_1
\Psi_1), \
\ \ \ \partial_z x^2 = \frac{1}{2} (\bar{\Phi}_2\bar{\Psi}_2 -
\Phi_1 \Psi_1),
\\
\partial_z x^3 = \frac{1}{2} (\bar{\Phi}_2 \Psi_1 + \Phi_1
\bar{\Psi}_2),
\ \ \ \ \partial_z x^4 = \frac{i}{2} (\bar{\Phi}_2 \Psi_1 - \Phi_1
\bar{\Psi}_2).
\end{split}
\end{equation}

The periodicity of the quantities $x^k_z, k=1,2,3,4$, with respect to
$\Lambda$ only implies that
$$
\Psi_1(z+\gamma_i) = \kappa_i \Psi_1(z), \ \ \
\Psi_2(z+\gamma_i) = \bar{\kappa}_i \Psi_2(z), \ \ \
$$
$$
\Phi_1(z+\gamma_i) = \frac{1}{\kappa_i} \Phi_1(z), \ \ \
\Phi_2(z+\gamma_i) = \frac{1}{\bar{\kappa}_i} \Phi_i(z), \ \ \ i=1,2.
$$
However these vector functions satisfy differential equations with
periodic coefficients:
$$
\D \Psi = 0, \ \ \D^\vee \Phi = 0
$$
which implies that $\kappa_i = \bar{\kappa}_i, i=1,2$, i.e. the
multipliers are real-valued.
Therefore we have
$$
\Psi(z+\gamma_i) = \hat{\kappa}_i \Psi(z), \ \ \ \Phi(z+ \gamma_i)
= \frac{1}{\hat{\kappa}_i} \Phi(z), \ \ \ i=1,2,
$$
and we may apply the results from the previous section.

Let a deformation of $L$, $L^\vee$ be a
sum of the following infinitesimal Darboux transformations:
$$
\partial_\tau = \partial_{\tau_1} + \partial_{\tau_2},
$$
where $\partial_{\tau_1} $  is generated by the following pair of
solutions:
$$
\Psi^D=\left( \begin{array}{c}  \Psi_1  \\ \Psi_2  \end{array}
\right), \ \ \Phi^D=\left( \begin{array}{c}  \bar\Phi_2  \\
-\bar\Phi_1
\end{array} \right)
$$
and  $\partial_{\tau_2} $  is generated by:
$$
\Psi^D=\left( \begin{array}{c}  \bar\Psi_2  \\ -\bar\Psi_1
\end{array} \right), \ \ \Phi^D=\left( \begin{array}{c}  \Phi_1
\\  \Phi_2
\end{array} \right).
$$

\begin{theorem}
\label{theorem2}
1. The flow $\partial_\tau$ is isospectral and respects the reality
conditions.  Moreover
$$
\partial_\tau U= \Phi_1\bar\Psi_1-\bar\Phi_2\Psi_2, \ \
\partial_\tau \bar U= \bar\Phi_1\Psi_1-\Phi_2\bar\Psi_2, \ \
$$

2. Let us assume, that the kernels $\omega$ and $\omega^\vee$ for the
functions $\Psi$ and $\Phi$,
are normalized by: $\omega(0)=0$, $\omega^\vee(0)=0$. Then
\begin{equation}
\begin{split}
\partial_\tau \Psi_1=(x^3-ix^4)\Psi_1 - i (x^1-ix^2) \bar\Psi_2  \\
\partial_\tau \Psi_2=(x^3-ix^4)\Psi_2  + i (x^1-ix^2) \bar\Psi_1 \\
\partial_\tau \Phi_1=(x^3+ix^4)\Phi_1 - i (x^1-ix^2) \bar\Phi_2  \\
\partial_\tau \Phi_2=(x^3+ix^4)\Phi_2  + i (x^1-ix^2) \bar\Phi_1.
\end{split}
\end{equation}
\end{theorem}

{\it Proof}. The first part follows directly from
Theorem~\ref{lemma-iso}.

For the derivatives of the coordinate functions  we have
\begin{equation}
\begin{split}
-i(x^1-ix^2)_z=  \Phi_1\Psi_1 \\
-i(x^1+ix^2)_z=   \bar\Phi_2\bar\Psi_2 \\
(x^3+ix^4)_z=  \Phi_1\bar\Psi_2 \\
(x^3- ix^4)_z=   \bar\Phi_2\Psi_1,
\end{split}
\end{equation}
and for the flow $\partial_{\tau_1}$
$$
\partial_z\omega= \bar\Phi_2 \Psi_1, \ \  \partial_z\omega^\vee =
\Phi_1 \Psi_1.
$$
Therefore
$$
\omega=(x^3- ix^4), \ \  \omega^\vee = -i(x^1-ix^2),
$$
\begin{equation}
\begin{split}
\partial_{\tau_1} \Psi_1=(x^3-ix^4)\Psi^D_1  = (x^3-ix^4)\Psi_1 \\
\partial_{\tau_1} \Psi_2=(x^3-ix^4)\Psi^D_2  = (x^3-ix^4)\Psi_2  \\
\partial_{\tau_1} \Phi_1=-  i (x^1-ix^2) \Phi^D_1  =  - i (x^1-ix^2)
\bar\Phi_1  \\
\partial_{\tau_1} \Phi_2= - i (x^1-ix^2) \Phi^D_2  =   i (x^1-ix^2)
\bar\Phi_2 .
\end{split}
\end{equation}
For the flow $\partial_{\tau_2}$
$$
\partial_z\omega= \Phi_1 \Psi_1, \ \  \partial_z\omega^\vee =  \Phi_1
\bar\Psi_2.
$$
Therefore
$$
\omega=-i(x^1- ix^2), \ \  \omega^\vee = (x^3+ix^4),
$$
\begin{equation}
\begin{split}
\partial_{\tau_2} \Psi_1=-i(x^1-ix^2)\Psi^D_1 = - i (x^1-ix^2)
\bar\Psi_2  \\
\partial_{\tau_2} \Psi_2=-i(x^1-ix^2)\Psi^D_2 =  i (x^1-ix^2)
\bar\Psi_1 \\
\partial_{\tau_2} \Phi_1=(x^3+ix^4)\Phi^D_1 =(x^3+ix^4)\Phi_1 \\
\partial_{\tau_2} \Phi_2=(x^3+ix^4)\Phi^D_2=(x^3+ix^4)\Phi_2 .
\end{split}
\end{equation}

\begin{theorem}
\label{theorem3} This action generate the following infinitesimal
conformal transformations of the immersed surface:
\begin{equation}
\label{generator}
\begin{split}
\partial_{\tau}x^1 = 2 x^1 x^3  \\
\partial_{\tau}x^2 = 2 x^2 x^3  \\
\partial_{\tau}x^3 = (x^3)^2 -(x^1)^2 - (x^2)^2 - (x^4)^2 \\
\partial_{\tau}x^4 = 2 x^4 x^3.
\end{split}
\end{equation}
\end{theorem}

{\it Proof}. It is sufficient to check that
\begin{equation}
\begin{split}
\partial_{\tau}\partial_z x^1 = 2 [\partial_z x^1] x^3 + 2 x^1
[\partial_z x^3]   \\
\partial_{\tau}\partial_z x^2 = 2 [\partial_z x^2] x^3 + 2 x^2
[\partial_z x^3] \\
\partial_{\tau}\partial_z x^3 = 2 x^3 [\partial_z x^3] - x^1
[\partial_z x^1]
 - x^2 [\partial_z x^2] - x^4 [\partial_z x^4] \\
\partial_{\tau}\partial_z x^4 = 2 [\partial_z x^4] x^3 + 2 x^4
[\partial_z x^3].
\end{split}
\end{equation}

We obtain by straightforward computations that
$$
\partial_\tau\partial_z x^1 = \frac{i}{2} (
[\partial_\tau\bar{\Phi}_2]\bar{\Psi}_2 +
\bar{\Phi}_2[\partial_\tau\bar{\Psi}_2] +
[\partial_\tau\Phi_1]\Psi_1 + \Phi_1 [\partial_\tau\Psi_1] )=
$$
$$
=\frac{i}{2}([(x^3-ix^4)\bar\Phi_2  - i (x^1+ix^2)\Phi_1]
\bar{\Psi}_2+ \bar{\Phi}_2 [(x^3+ix^4)\bar\Psi_2  - i
(x^1+ix^2)\Psi_1] +
$$
$$
+[(x^3+ix^4)\Phi_1 - i (x^1-ix^2) \bar\Phi_2 ]\Psi_1 +
\Phi_1[(x^3-ix^4)\Psi_1 - i (x^1-ix^2)\bar\Psi_2 ])=
$$
$$
=\frac{i}{2}( 2 x^3 ( \bar\Phi_2\bar{\Psi}_2 + \Phi_1 \Psi_1 ) + 2
x^1 (-i \Phi_1 \bar{\Psi}_2 -i \bar\Phi_2 \Psi_1) ) = 2 x^3
\partial_z x^1 + 2 x^1 \partial_z x^3,
$$
$$
\partial_\tau\partial_z x^2 = \frac{1}{2} (
[\partial_\tau\bar{\Phi}_2]\bar{\Psi}_2 +
\bar{\Phi}_2[\partial_\tau\bar{\Psi}_2] -
[\partial_\tau\Phi_1]\Psi_1 - \Phi_1 [\partial_\tau\Psi_1] )=
$$
$$
=\frac{1}{2}([(x^3-ix^4)\bar\Phi_2  - i (x^1+ix^2)\Phi_1]
\bar{\Psi}_2+ \bar{\Phi}_2 [(x^3+ix^4)\bar\Psi_2  - i
(x^1+ix^2)\Psi_1] -
$$
$$
-[(x^3+ix^4)\Phi_1 - i (x^1-ix^2) \bar\Phi_2 ]\Psi_1 -
\Phi_1[(x^3-ix^4)\Psi_1 - i (x^1-ix^2)\bar\Psi_2 ])=
$$
$$
=\frac{1}{2}( 2 x^3 ( \bar\Phi_2\bar{\Psi}_2 - \Phi_1 \Psi_1 ) + 2
x^2 (\Phi_1 \bar{\Psi}_2 + \bar\Phi_2 \Psi_1) ) = 2 x^3 \partial_z
x^2 + 2 x^2 \partial_z x^3,
$$
$$
\partial_\tau\partial_z x^3 = \frac{1}{2} (
[\partial_\tau\bar{\Phi}_2]{\Psi}_1 +
\bar{\Phi}_2[\partial_\tau{\Psi}_1] +
[\partial_\tau\Phi_1]\bar\Psi_2 + \Phi_1 [\partial_\tau\bar\Psi_2]
)=
$$
$$
=\frac{1}{2}([(x^3-ix^4)\bar\Phi_2  - i (x^1+ix^2)\Phi_1]
{\Psi}_1+ \bar{\Phi}_2 [(x^3-ix^4)\Psi_1 - i (x^1-ix^2)\bar\Psi_2
] +
$$
$$
+[(x^3+ix^4)\Phi_1 - i (x^1-ix^2) \bar\Phi_2 ]\bar\Psi_2 +
\Phi_1[(x^3+ix^4)\bar\Psi_2  - i (x^1+ix^2)\Psi_1])=
$$
$$
=-\frac{i}{2}( 2 x^1 ( \bar\Phi_2\bar{\Psi}_2 + \Phi_1 \Psi_1 ) )-
\frac{1}{2}( 2 x^2 ( \bar\Phi_2\bar{\Psi}_2 - \Phi_1 \Psi_1 ) ) +
$$
$$
+ \frac{1}{2}( 2 x^3 ( \bar\Phi_2 \Psi_1 +  \Phi_1  \bar\Psi_2 ))
- \frac{i}{2}( 2 x^4 (\bar\Phi_2 \Psi_1 - \Phi_1  \bar\Psi_2 ) )=
$$
$$
= - 2 x^1 \partial_z x^1 - 2 x^2 \partial_z x^2  + 2 x^3
\partial_z x^3 - 2 x^4 \partial_z x^4
$$
and
$$
\partial_\tau\partial_z x^4 = \frac{i}{2} (
[\partial_\tau\bar{\Phi}_2]{\Psi}_1 +
\bar{\Phi}_2[\partial_\tau{\Psi}_1] -
[\partial_\tau\Phi_1]\bar\Psi_2 - \Phi_1 [\partial_\tau\bar\Psi_2]
)=
$$
$$
=\frac{i}{2}([(x^3-ix^4)\bar\Phi_2  - i (x^1+ix^2)\Phi_1]
{\Psi}_1+ \bar{\Phi}_2 [(x^3-ix^4)\Psi_1 - i (x^1-ix^2)\bar\Psi_2
] -
$$
$$
-[(x^3+ix^4)\Phi_1 - i (x^1-ix^2) \bar\Phi_2 ]\bar\Psi_2 -
\Phi_1[(x^3+ix^4)\bar\Psi_2  - i (x^1+ix^2)\Psi_1])=
$$
$$
= \frac{i}{2}( 2 x^3 ( \bar\Phi_2 \Psi_1 - \Phi_1  \bar\Psi_2 )) +
\frac{1}{2}( 2 x^4 (\bar\Phi_2 \Psi_1 + \Phi_1  \bar\Psi_2 ) ) = 2
x^3 \partial_z x^4 + 2 x^4 \partial_z x^3.
$$
This completes the proof of Theorem \ref{theorem3}.

It is known that the group of conformal transformations of $\bar{\R}^4$
is generated by rotations,
dilations, translations and inversions. The spectral curve is evidently
preserved by rotations,
translations and dilations. Any inversion is conjugate by these
transformations to the inversion
with the generator of the form (\ref{generator}). It is also clear that
the set of conformal transformations
which map a given torus into a compact surface is path-connected.
Therefore
Theorems~\ref{theorem2} and \ref{theorem3} imply

\begin{theorem}
\label{theorem4}
Any conformal transformation of $\bar{\R}^4$ which maps a torus $T
\subset \R^4$ into
a compact surface preserves the multiplier set of the corresponding
Dirac operator.
\end{theorem}

As we mentioned in \S \ref{sec2} there is some freedom in choosing
the Dirac operators for a given torus $T \subset \R^4$. The
possible potentials form a discrete set $P_T$. Conformal
transformations of $\bar{\R}^4$ which map the torus $T$ into a
compact surface form a path-connected component $G_T$ of the
conformal group $SO(5,1)$ and the action of $G_T$ on $P_T$ is
defined. Therewith Theorem~\ref{theorem4} reads that this action
and, in particular, any continuous family $f_t \subset G_T$
preserves the multiplier sets.

\section{Examples}
\label{sec5}

\subsection{Spectral curves of the Clifford tori}
\label{subsec5.1}

Following \cite{T3} we present a pair of tori in $\R^4$ which are
related by a conformal transformation,
have the same multiplier set however their spectral curves are
different.

In fact these are the Clifford tori in $S^3$ and $\R^3$ which are
related by a stereographical projection.
the famous Willmore conjecture reads that for tori in $\R^3$ the
Willmore functional attains its minimum
on the Clifford torus and its images under conformal transformations.
This conjecture is generalized for tori in $\R^4$
and reads actually the same: the Willmore functional
$$
{\cal W}(M) = \int_M |{\bf H}|^2 d\mu
$$
attains its minimum on the Clifford torus in
$S^3$ and its conformal images (here ${\bf H}$ is the mean curvature
vector).

A) {\it The Clifford torus in the unit three-sphere} $S^3 \subset
\R^4$
is defined by the equations
$$
\left(x^1\right)^2 + \left(x^2\right)^2 = \frac{1}{2}, \ \ \
\left(x^3\right)^2 + \left(x^4\right)^2 = \frac{1}{2}
$$
where
the sphere is defined by the equation
$$
\left(x^1\right)^2 + \left(x^2\right)^2 +
\left(x^3\right)^2 + \left(x^4\right)^2 = 1.
$$

Its spectral curve is $\Gamma$ is the sphere $\C P^1 = \bar{\C}$
with two marked points $\lambda=0,\infty$ added at the
compactification.

The potential of the Dirac operator $\D$ equals
$$
U = \frac{1+i}{4}.
$$
The zero Floquet--Bloch functions are parameterized by the points of
$\Gamma^\prime = \Gamma \setminus\{\lambda=0,\infty\}$
and are glued into the Baker--Akhiezer function
$\psi$ which is is meromorphic on $\Gamma^\prime$ and
has at the marked points (``infinities'') the following asymptotics:
$$
\psi \approx
\left(
\begin{array}{c}
e^{\lambda  z} \\ 0
\end{array}
\right) \ \mbox{as $\lambda \to \infty$};
\ \
\psi \approx
\left(
\begin{array}{c}
0 \\ e^{- \frac{|u|^2}{\lambda} \bar{z}}
\end{array}
\right) \ \mbox{as $\lambda \to 0$}
$$
where $u=\frac{1+i}{4}$.

The Clifford torus is defined via the Weierstrass representation
by the functions
$$
\psi_1 = \psi_2 = \frac{1}{\sqrt{2}}
\exp \left( -\frac{i(x+y)}{2} \right), \ \
\varphi_1 = -\varphi_2 = -\frac{1}{2\sqrt{2}}
\exp \left( \frac{i(y-x)}{2} \right).
$$

Although the Dirac operator for a torus in $\R^4$ is defined up to
transformations (\ref{curvetransform}) and replacing $\D$ by
$\D^\vee$ it is easy to check that for all such choices the spectral
curve of the Clifford torus stays smooth and coincide with $\C P^1$.

B) {\it The Clifford torus in $\R^3$} is the image of the Clifford
torus in
$S^3 \subset \R^4$ under a stereographic projection
$$
(x^1,x^2,x^3,x^4) \to \left(\frac{x^1}{1-x^4},\frac{x^2}{1-x^4},
\frac{x^3}{1-x^4}\right), \ \ \ \sum_k (x^k)^2 = 1,
$$
which is extended to a conformal transformation of $\bar{\R}^4$.

Its spectral curve is $\Gamma$ is the sphere $\C P^1 = \bar{\C}$
with two marked points $\lambda=0,\infty$ added at the
compactification and two double points obtained by gluing together
the points from the following pairs:
$$
\left(\frac{1+i}{4},\frac{-1+i}{4}\right)
\ \ \ \mbox{and} \ \ \
\left(-\frac{1+i}{4},\frac{1-i}{4}\right).
$$
The potential of the Dirac operator $\D$ equals
$$
U = \frac{\sin y}{2\sqrt{2}(\sin y  - \sqrt{2})}
$$
The zero Floquet--Bloch functions are parameterized by the points
of $\Gamma^\prime = \Gamma \setminus\{\lambda=0,\infty\}$ and are
glued into the Baker--Akhiezer function which is uniquely defined
by the following conditions:

1) it is meromorphic on $\Gamma^\prime$ and has the following
asymptotics
$$
\psi \approx
\left(
\begin{array}{c}
e^{\lambda  z} \\ 0
\end{array}
\right) \ \mbox{as $\lambda \to \infty$};
\ \
\psi \approx
\left(
\begin{array}{c}
0 \\ e^{- \frac{|u|^2}{\lambda} \bar{z}}
\end{array}
\right) \ \mbox{as $\lambda \to 0$}
$$
where $u=\frac{1+i}{4}$;

2) it has three poles
$\Gamma^\prime$ which are
independent on $z$ and have the form
$$
p_1 = \frac{-1+i + \sqrt{-2i-4}}{4\sqrt{2}}, \ \
p_2 = \frac{-1+i - \sqrt{-2i-4}}{4\sqrt{2}}, \ \
p_3 = \frac{1}{\sqrt{8}}.
$$

The Clifford torus in $\R^3$ is constructed via the Weierstrass
representation
from the function
$$
\psi = \psi\left(z,\bar{z},\frac{1-i}{4}\right)
$$
which has the form
$$
\psi_1(z,\bar{z},\lambda) = e^{\lambda z
-\frac{|u|^2}{\lambda}\bar{z}}
\left(
q_1 \frac{\lambda}{\lambda-p_1} + q_2 \frac{\lambda}{\lambda-p_2} +
(1-q_1-q_2)\frac{\lambda}{\lambda-p_3}\right),
$$
$$
\psi_2(z,\bar{z},\lambda) = e^{\lambda z
-\frac{|u|^2}{\lambda}\bar{z}}
\left(
t_1 \frac{p_1}{p_1-\lambda} + t_2
\frac{p_2}{p_2-\lambda} +
(1-t_1-t_2)\frac{p_3}{p_3-\lambda}\right)
$$
where $u=\frac{1+i}{4}$. The functions $q_1,q_2,t_1,t_2$ depend
only on $y$ and $2\pi$-periodic with respect to $y$ and found from
the following conditions
$$
\psi\left(z,\bar{z},\frac{1+i}{4}\right) =
\psi\left(z,\bar{z},\frac{-1+i}{4}\right), \ \ \
\psi\left(z,\bar{z},-\frac{1+i}{4}\right) =
\psi\left(z,\bar{z},\frac{1-i}{4}\right).
$$

We see that the spectral curve of the Clifford torus in $S^3$ is
smooth and the spectral curve of the Clifford torus in $\R^3$ has
a pair of double points. Their multiplier sets are the same.
Therefore the spectral curves can be deformed by conformal
transformations however, by Theorem 1, such a deformation may only
consist in gluing together or ungluing multiple points with the
multipliers $\kappa_i/\hat{\kappa}_i  = 1$ for $\Gamma$ and
$\kappa_i \hat{\kappa}_i  = 1$ for $\Gamma^\vee$.

In particular, we see that in a completely conformal setting
adopted, for instance, in \cite{FLPP} it is impossible to define
the spectral curve $\Gamma$ as we did above and its definition
always needs the analytical theory of differential operators.

Let us point out, that in the both cases discussed above the potential
does not depend on the variable $x$: $U(x,y)=U(y)$. For such potentials
the Floquet solution has the following form:
$$
\psi(\lambda,x,y)=\left(\begin{array}{c} \tilde\psi_1(\lambda,y) \\
\tilde\psi_2(\lambda,y)\end{array}\right)e^{kx}, \ \ \
\varphi(\lambda,x,y)=\left(\begin{array}{c} \tilde\psi_2(\sigma\lambda,y) \\
\tilde\psi_1(\sigma\lambda,y)\end{array}\right)e^{-kx},
$$
where $\tilde\psi(\lambda,y)$ denotes the Floquet solution for the
1-dimensional Dirac operator
$$
\tilde L \tilde\psi(\lambda,y)=-k\tilde\psi(\lambda,y),\ \
\tilde L =\left[\begin{array}{cc} i\partial_y & -2\bar U \\
-2U &  -i\partial_y \end{array}\right], \ \
\tilde\psi(\lambda,y)=\left(\begin{array}{c} \tilde\psi_1(\lambda,y) \\
\tilde\psi_2(\lambda,y)\end{array}\right).
$$
$\tilde L$ is the auxiliary operator for the self-focusing nonlinear
Schr\"odinger equation (NLS). The spectral curve $\Gamma$ for
$\tilde L$ is a two-sheeted covering of the $k$-plane.  All
branching points of $ \Gamma$ lie outside the real line and form
complex conjugate pairs. Here $\sigma$ denotes the transposition of
sheets.

The kernels $\omega(\lambda,z)$, $\omega^\vee(\mu,z)$ can be easily
calculated explicitly. Assuming
$$
\Psi^D(z)=e^{k_1x} \tilde\Psi^D(y), \ \ \Phi^D(z)=e^{-k_1x} \tilde\Phi^D(y),
$$
and integrating over $x$ we obtain

$$
\omega(\lambda,z) =\frac{e^{(k(\lambda)-k_1)x}}{k(\lambda)-k_1}
\left( \tilde\Phi^D_1(y)\tilde\psi_1(\lambda,y)-
\tilde\Phi^D_2(y)\tilde\psi_2(\lambda,y)\right)
$$
$$
\omega^\vee(\mu,z) =\frac{e^{(k_1-k(\mu))x}}{k_1-k(\mu)}
\left(\tilde\psi_2(\sigma\mu,y)\tilde\Psi^D_1(y)-
\tilde\psi_1(\sigma\mu,y)\tilde\Psi^D_2(y)
\right)
$$

The spectral curve of the Clifford torus corresponds to the
so-called whiske\-red torus of the NLS equation \cite{McL}. It
means, that the Liouville torus is the product $S^1\times \hat T$,
where $\hat T$ denotes the one-point  compactification of the
product $S^1\times{\mathbb R}^1$. The Liouville torus is represented
as a compact set in the phase space. The $y$-dynamics of the NLS
equation and the phase gauging flow $U\rightarrow Ue^{i\phi}$ are
orthogonal to the infinite direction, and the time-evolution
corresponds to the motion along the infinite cycle. The limiting
point corresponds to the Clifford torus in $S^3$. The conformal flow
is proportional to the time evolution plus phase gauging, but the
coefficient between these flows became infinite near the fixed
point, therefore using the conformal transformations one reaches the
fixed point at finite time.

The whiskered tori are very important in the theory of the NLS
equation, because they are the principal source of instability with
respect to small perturbations, and, as a corollary, generate
numerical chaos \cite{AH1,AH2,McL}.

\subsection{The spectral curve of the one-dimensional
Schr\"odinger operator with a complex potential}
\label{subsec5.2}

We consider the one-dimensional Schr\"odinger operator
$$
{\cal L} = -\frac{d^2}{dx^2} + u(x)
$$
associated with the KdV equation. We assume that $u(x)$ is a real
periodic finite-gap potential: $u(x+T)=u(x)$. Let us consider the
$g$-gap potential. The spectral gaps are $(-\infty,E_0)$,
$(E_1,E_2)$, \ldots, $(E_{2g-1},E_{2g})$. The surface $\Gamma$
parameterizing the Floquet--Bloch solutions is the two-sheeted
covering of the complex plane ramified at the branch points
$\infty$, $E_k$, $k=0,\ldots,2g$. Compactifying $\Gamma$ by adding
an infinite point one obtains a smooth algebraic curve of genus
$g$. The Floquet--Bloch function $\psi(x,P), P in \Gamma$, is
defined as a function on the Riemann surface $\Gamma$:
\footnote{In the terms of \S \ref{sec2} this complex curve may be
considered as the complex curve of the two-dimensional
Schr\"odinger operator $\partial_y - \partial^2_x + u(x)$ on the
zero energy level.}
$$
w^2 = Q(E) = (E-E_0)\dots(E-E_{2g}),
$$
with the asymptotic
$$
\psi(P,x) = e^{i\sqrt{E}x}\left(1 +
O\left(\frac{1}{\sqrt{E}}\right)\right).
$$
It is assumed that it is meromorphic for finite $E$ and it is
uniquely fixed by its $g$ poles lying at the points
$(\gamma_j,\sqrt{Q(\gamma_j}), j=1,\dots,g$.

The operator ${\cal L}$ has infinitely many resonant points
$\tilde E_j\in \C$, namely the energy levels at which  the
monodromy operator becomes equal to $\pm I$, where $I$ is the unit
$2\times2$ matrix. These points can be treated as degenerate gaps
of zero length. Denote the preimages of $\tilde E_j$ in $\Gamma$
by $\tilde E_j^{\pm}$. The pairs  $\tilde E_j^{\pm}$ are called
resonant pairs.

The immersion of $\Gamma$ into $\C^2$ is defined by
$\M(\gamma)=(E,\kappa(\gamma))$, where $\kappa$ is the multiplier
with respect to the shift $x\rightarrow x+T$. The curve
$\M(\Gamma)$ has infinitely many self-intersections:
$\M(E_j^+)=\M(E_j^-)$ for all $j$.

The set of all real potentials corresponding to this curve
$\M(\Gamma)$ is a real $g$-dimensional torus, and the spectral
curve parameterizing the variety of Floquet-Bloch functions
coincides with $\Gamma$ for all members of this family. But if we
consider the complex potentials corresponding to the curve
$\M(\Gamma)$, the situation drastically changes. These potentials
form an infinite-parametric family, and for generic members of
this family the Floquet-Bloch solutions are parameterized by a
curve, obtained from $\Gamma$ by gluing together all resonant
pairs. From the algebro-geometrical point of view these potentials
are infinite-gap. If only $k$ resonant pairs are glued, the
corresponding potentials form a family of complex dimension $g+k$,
and the arithmetic genus of the parameterizing spectral curve is
equal to $g+k$.

At the analytic level the procedure of gluing of a resonant pair
means, that an extra pole is added to the divisor, and
simultaneously an extra linear relation is imposed on the wave
function $\psi(\gamma,x)$: $\psi(\tilde E_j^+,x)=\psi(\tilde
E_j^-,x)$ for all $x$.

The opposite operation of ungluing the pair $(\tilde E_j^+$,
$\tilde E_j^-)$ corresponds to the following degeneration: one of
the divisor points $\gamma_l$ tends to either $\tilde E_j^+$ or
$\tilde E_j^-$. It is easy to check, then the residue at
$\gamma_l$ is proportional to $\gamma_l-\tilde E_j$,  and after
taking the limit the pole vanishes, but the values of the wave
functions at the points $\tilde E_j^+$ and $\tilde E_j^-$ become
different.

\end{document}